[1994/12/01]
\documentclass{amsart}
\chardef\bslash=`\\ 

\usepackage[square,sort,authoryear, colon]{natbib}
\bibpunct{[}{]}{;}{a}{}{,}

\usepackage{amsthm,thmtools}

\usepackage{amssymb}

\usepackage{color}
\usepackage{graphicx}

\usepackage{caption}
\usepackage{wrapfig}
\usepackage{lscape}
\usepackage{rotating}

\usepackage{overpic}

\usepackage{adjustbox}

\usepackage{hvfloat,lipsum}

\usepackage{caption}
\usepackage{subcaption}

\usepackage{appendix}

\usepackage{url}

\usepackage{accents}

\newtheorem{theorem}{\bf Theorem}[section]
\newtheorem{proposition}[theorem]{\bf Proposition}

\newtheorem{conjecture}[theorem]{\bf Conjecture}
\newtheorem{example}[theorem]{\bf Example}

\newenvironment{myquote}[1]%
  {\list{}{\leftmargin=#1\rightmargin=#1}\item[]}%
  {\endlist}

\usepackage{hyperref}
\usepackage{csquotes}

\theoremstyle{definition}
\newtheorem*{definition}{Definition}
\newtheorem*{remark}{Remark}

\definecolor{gray}{rgb}{0.8,0.8,0.8}

\newcommand{\obsgray}[1]{\textcolor{gray}{#1}}

\newcommand{\sieved}[1]{\textcolor{gray}{#1}}

\newcommand{\e}{\mathrm{e}}

\newcommand{\propositionref}[1]{Proposition~\ref{#1}}

\newcommand{\figref}[1]{Figure~\ref{#1}}

\newcommand{\tabref}[1]{Table~\ref{#1}}

\usepackage{mathtools}

\DeclareMathOperator{\li}{li}

\makeatletter
\def\imod#1{\allowbreak\mkern10mu({\operator@font mod}\,\,#1)}
\makeatother

\newcommand{\eval}[2][\right]{\relax
  \ifx#1\right\relax \left.\fi#2#1\rvert}

\title[A first-principles perspective on the primes]{Defining the prime numbers prior to the integers: A first-principles approach to the distribution of primes}

\author{Kolbj{\o}rn Tunstr{\o}m}
\email{kolbjorn@chalmers.se}
\address{Department of Physics, Chalmers University of Technology, 41296 Gothenburg, Sweden}

\begin{document}
\begin{abstract}
While the prime numbers have been subject to mathematical inquiry since the ancient Greeks, the accumulated effort of understanding these numbers has---as Marcus du Sautoy recently phrased it---'not revealed the origins of what makes the primes tick.' Here, we suggest that a resolution to this long-standing conundrum is attainable by defining the primes \emph{prior} to the natural numbers---as opposed to the standard number theoretical definition of primes where these numbers derive \emph{from} the natural numbers. 
The result is a first-principles perspective on the primes that exposes and explains the 'origins' of their distribution and their mathematical properties and provides an intuitive as well as pedagogical approach to the primes with the potential to impact our thinking about these age-old mathematical objects. A few immediate outcomes of this perspective are another proof of the fundamental theorem of arithmetic, a probabilistic model of primes sharing as well as explaining their subrandom correlation structure, and an equivalent formulation of the Riemann hypothesis. 
\end{abstract}

 \maketitle


\section{Introduction}
As the building blocks of the natural numbers, the prime numbers rank among the most fundamental objects in the whole of mathematics. Remarkably, their origin as mathematical objects in their own right 
dates back more than two thousand years---the earliest surviving mathematical treatment of primes being found in Euclid's Elements (c. 300 BCE). The modern study of primes, on the other hand, and in particular how the primes  are distributed in the natural numbers, began its 
development first in the 17th century and onwards, eventually evolving into the branch of mathematics now recognised as analytic number theory. This development was especially influenced by Riemann's legendary foray into prime numbers---the 1859 manuscript \emph{\"Uber die Anzahl der Primzahlen unter einer gegebenen Gr\"osse}---which,  
only nine pages thin, outlined a set of monumental ideas that now underly much of the established theory on primes, including landmarks such as the prime number theorem. The quest for a deeper grasp of the primes' mathematical nature is still ongoing, and while obviously a mature field of mathematical research, the theory on primes has progressed profoundly even in recent times, as exemplified by renowned contemporary results such as the Green-Tao theorem on arithmetic progressions in the primes~\citep{Green:2008}, Zhang's theorem on bounded gaps between primes~\citep{Zhang:2014}, Helfgott's proof of the weak Goldbach conjecture~\citep{Helfgott:2013}, and a recent conjectural result on unexpected biases in the distribution of consecutive primes~\citep{Oliver:2016}. 

Despite the extensive body of theory and the sophisticated mathematical technologies developed to study and prove various properties of the primes---not to mention the time span of their study---the perhaps most fundamental question one could ask about the primes is still lingering unanswered, namely, what 
underlies 
the distribution of prime numbers? The backdrop of this question is the status quo that no satisfactory explanation exists for what is often portrayed as a mystery \citep{QuotesWeb, Luque:2009, Sautoy:2017}: 
The sequence of primes is deterministic, but the primes appear to be scattered almost randomly throughout the natural numbers---so much so that that the most accurate models of primes are random models~\citep{Tao:web:RandomModels}. In aggregate, though, the randomness translates into precise regularity, since the  prime number theorem guarantees that the asymptotic density of primes below $x$ equals $1/\log x$. To see this long-standing conundrum resolved, one would ultimately want an explanation rooted in first principles. As we establish in this paper, a slight shift of perspective on the primes is enough to accomplish this milestone. 

The initial motivation of an altered perspective  stems from a simple reflection about the primes: While the prime numbers are renowned as the multiplicative building blocks of the natural numbers, they are normally defined deconstructively in terms of what they build: 
%
%
%
%
\emph{A prime is a natural number greater than 1 whose only positive divisors are 1 and itself}. But one could reasonably argue that a more logical arrangement would be the building blocks appearing prior to the whole---not the other way around. In view of this, the obvious question emerging is whether there are alternative definitions that allow the primes to precede the natural numbers? And if so, would that bring about any new insights? In this inquiry, we affirmatively answer both questions.
%
%
%
%
%
%

To begin with, one definition that fits our ambition is due to Euclid (\textit{Elements: Book VII: Definition 11}) and reads: 
\begin{definition}[Euclid] 
    A prime number is that which is measured by a unit alone. 
\end{definition}

In the first part of this paper we exploit Euclid's definition to formulate a \emph{constructive} prime number generator, where the primes arise prior to the natural numbers. 
An immediate implication of this first-principles approach is that the
distribution of primes is easily understood in terms of an ever-increasing combination of recursively defined prime periodic sequences, 
%
where the initial generating element is a periodic sequence with period one unit, naturally represented by the 
sequence $1,1,1,\dots$. The natural numbers are defined subsequently in terms of prime periodic sequences, and the reversed order of appearance results in a novel elementary proof of the fundamental theorem of arithmetic. 
%
%
It is worth mentioning that this perspective ties together addition and multiplication naturally in the sense that the unit element 1 can be viewed as the underlying generator for the natural numbers in both cases.

The recursive nature of the primes reveals a well-defined, discrete, structure of the distribution of primes in relation to the natural numbers, namely that the distribution of primes between the $k$th and $k+1$st primes squared is fully determined by the $k$ first prime periodic sequences. Obviously, this discrete structure---which yields a complete subdivision of the natural numbers---is apparent also from the sieve of Eratosthenes.
The latter part of the paper is concerned both with how the reputed randomlike distribution of the primes results from the recursive build-up of the stated structure, as well as with explaining in what sense the primes are randomly distributed. In fact, by defining and analysing a random model of the primes that accounts for their structure, we learn that the sequence of primes belongs to a category of sequences known as subrandom, recognised by having variance less than that of a corresponding random process. Moreover, the model accurately predicts as well as explains the observed subrandomness of the primes, which entirely derives from their recursive structure. 
%
%
%
%
%
%
%
%
%
%
%
Consequently, our random model seems
a valuable
tool for building bottom-up understanding of why many conjectures about primes should be true in the first place, including the Riemann Hypothesis and the Hardy-Littlewod k-tuple conjecture, as well as the conjectured results in~\citep{Oliver:2016}.



The aspiration of this paper is to bring to light a first-principles perspective on the distribution of prime numbers that could contribute to our thinking about these age-old mathematical objects. Having said that, we are well aware of the sentiment expressed almost 100 years ago by Hardy and Littlewood, and still alive today, that "...in pure mathematics, and in The Theory of Numbers in particular, 'it is only proof that counts'" \citep[p. 68]{hardy1923}. While both propositions and proofs are included here for relevant results, and as such live up to this sentiment, we add the caveat that our emphasis is nevertheless on lifting a different perspective on the primes, and not to prove some of the outstanding conjectures in number theory. As Tao argues in "What is good mathematics?" \citep{Tao:2007}, there are many aspects to mathematical quality, and it is our hope that several of these are expressed well enough in our work to admit thoughtful consideration. 
%
%
%

Due to its conceptual leaning, the paper is suitable for---and deliberately aimed at---a broad mathematical audience. For this reason, we have placed the weight on explanation, grounded in theory as well as in heuristics based on numerical evidence. In addition, to enforce the proposed bottom-up perspective, we sometimes use notation that number theorists might prefer to replace with established notation. Our presentation thus deviates from the bulk of contemporary mathematical expositions and should be approached accordingly. Regardless, our inquiry suggests that beyond a complete first-principles explanation of the distribution of primes---which closes a foundational gap in our understanding of these entities---the insights offered have the quality and potential to influence new proof strategies for problems related to the distribution of primes, presumably of serious interest to the mathematical community.
\section{A constructive prime generator}
\label{sec:constructivesieve}
In order to employ Euclid's definition to generate the prime numbers, and in turn the natural numbers, we first need to interpret \emph{what} is to be measured. Recognising that Euclid's definition is geometric in nature, we shall consider a geometrical perspective alongside a numerical one. Starting off with the geometrical point of view, let us assume a \emph{unit ruler} of arbitrary length, as illustrated uppermost in Figure~1A. 
\begin{figure}[t]
\centering 
\includegraphics[width=0.7\textwidth]{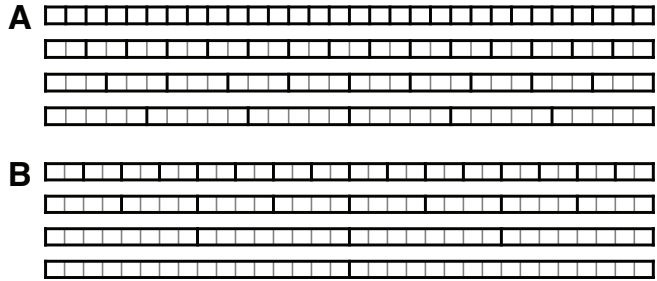}
\caption{\small
Geometrical version of the constructive prime generator. A: The set of prime rulers generated by the three initial steps of the first part of the constructive prime generator. The unit ruler is shown on top. Below are the prime rulers with units $p_1=2$, $p_2=3$, and $p_3=5$. B: The set of renormalised prime rulers generated by the three initial steps of the second part of the constructive prime generator, applied to the prime ruler with unit $p_1$, shown on top. Below are the renormalised prime rulers with units $p_1^2$, $p_1^3$, and $p_1^4$. 
}
\label{fig:rulers}
\end{figure}
Numerically, the unit ruler is paralleled by a sequence of 1s, which we denote the \emph{unit sequence}:
%
\begin{center}
\resizebox{0.7\columnwidth}{!}{%
\begin{tabular}{cccccccccccccccccc}
\obsgray{1}  & 1 & 1 & 1 & 1 & 1 & 1 & 1 & 1 & 1 & 1 & 1 & 1 & 1 & 1 & 1 & 1 & $\dots$
\end{tabular}%
}
\end{center}
The role of the unit ruler (or sequence) is to serve as a generating element for the set of prime numbers and eventually the natural numbers. 

Connecting back to Euclid's definition, what we want to measure is the distance from the beginning of the unit ruler to any position further down the ruler, in the process applying Euclid's definition to build the set of prime numbers and subsequently the natural numbers. We shall denote the resulting algorithm the \emph{constructive prime generator}. 
%
%
%
%
%
%
%
For the sake of clarity, we divide the constructive prime generator into two parts. The first part is responsible for constructing the set of primes, or more specifically, what we will call the set of \emph{prime rulers} (or numerically, \emph{prime sequences}), while the second part concerns constructing the set of \emph{renormalised prime rulers} (or numerically, \emph{renormalised prime sequences}). We will eventually see that any distance on the unit ruler larger than one unit can be measured in terms of rulers from these two sets---or equivalently, that any natural number greater than 1 is definable in terms of the combined set of prime sequences and renormalised prime sequences. 

The first part of the constructive prime generator is a simple iterative procedure, where each iteration step is executed as follows: We move along the unit ruler, beyond the first unit, until we reach a position measurable only in terms of the unit ruler, and not by any previously constructed \emph{prime rulers}. According to Euclid's definition, this distance defines a prime number, from which we construct a corresponding prime ruler. The first step of this procedure results in the initial prime ruler with unit $p_1$, which we then align with the unit ruler. Continuing, we obtain rulers with units $p_2$, $p_3$, and so on, as seen in Figure~1A. From the numerical perspective, we express the unfolding of the process as the repeated construction of \emph{prime sequences}; periodic sequences with repeating elements $1$ and $p_k$ and periods $p_k$,
$k\geq 1$, as illustrated in {\tabref{tab:table1}A}.
\begin{table}
\centering
\caption{\small \label{tab:table1}Numerical version of the constructive prime generator. A: The set of prime sequences generated by the first part of the constructive prime generator. The unit sequence is shown on top. Below are the prime sequences with periods $p_1=2$, $p_2=3$, $p_3=5$, and so on. B: The set of renormalised prime sequences generated by the second part of the constructive prime generator, applied to the prime sequence with period $p_1$, shown on top. Below are the renormalised prime sequences with periods $p_1^2$, $p_1^3$, $p_1^4$, and so on. 
}
\begin{tabular}{c}
\Large \textbf{\textsf{A}}
\end{tabular}
\quad	
\resizebox{0.85\columnwidth}{!}{%
\begin{tabular}{cccccccccccccccccc}
\obsgray{1} & 1 & 1 & 1 & 1 & 1 & 1 & 1 & 1 & 1 & 1 & 1 & 1 & 1 & 1 & 1 & 1 & $\dots$\\
\obsgray{$p_1$} & 1 & $p_1$ & 1 & $p_1$ & 1 & $p_1$ & 1 & $p_1$ & 1 & $p_1$ & 1 & $p_1$ & 1 & $p_1$ & 1 & $p_1$ & $\dots$\\
\obsgray{$p_2$} & 1 & 1 & $p_2$ & 1 & 1 & $p_2$ & 1 & 1 & $p_2$ & 1 & 1 & $p_2$ & 1 & 1 & $p_2$ & 1 & $\cdots$\\
\obsgray{$p_3$} & 1 & 1 & $1$ & 1 & $p_3$ & 1 & 1 & 1 & 1 & $p_3$ & 1 & 1 & 1 & 1 & $p_3$ & 1 & $\cdots$\\
\obsgray{\vdots}&\vdots&\vdots&\vdots&\vdots&\vdots&\vdots&\vdots&\vdots&\vdots&\vdots&\vdots&\vdots&\vdots&\vdots&\vdots&\vdots&$\cdots$ 
\end{tabular}%
}\\
\vspace{5mm}
\begin{tabular}{c}
{
\Large \textbf{\textsf{B}}
}
\end{tabular}
\quad	
\resizebox{0.85\columnwidth}{!}{%
\begin{tabular}{cccccccccccccccccc}
\obsgray{$p_1$} & 1 & $p_1$ & 1 & $p_1$ & 1 & $p_1$ & 1 & $p_1$ & 1 & $p_1$ & 1 & $p_1$ & 1 & $p_1$ & 1 & $p_1$ & $\dots$\\
\obsgray{$p_1$} & 1 & $1$ & 1 & $p_1$ & 1 & $1$ & 1 & $p_1$ & 1 & $1$ & 1 & $p_1$ & 1 & $1$ & 1 & $p_1$ & $\dots$\\
\obsgray{$p_1$} & 1 & $1$ & 1 & $1$ & 1 & $1$ & 1 & $p_1$ & 1 & $1$ & 1 & $1$ & 1 & $1$ & 1 & $p_1$ & $\dots$\\
\obsgray{$p_1$} & 1 & $1$ & 1 & $1$ & 1 & $1$ & 1 & $1$ & 1 & $1$ & 1 & $1$ & 1 & $1$ & 1 & $p_1$ & $\dots$\\
\obsgray{\vdots}&\vdots&\vdots&\vdots&\vdots&\vdots&\vdots&\vdots&\vdots&\vdots&\vdots&\vdots&\vdots&\vdots&\vdots&\vdots&\vdots&$\cdots$
\end{tabular}%
}
\end{table}

Evidently, the first part of the constructive prime generator straightforwardly generates the prime numbers in terms of the units of the prime rulers or the periods of the prime sequences---and is of course equivalent to the ancient and oft-used sieve of Eratosthenes. This means we now in principle have knowledge of what prime rulers are involved in measuring any distance on the unit ruler---or equivalently, what are the distinct prime factors of any natural number. But this only takes us half the way to the natural numbers, as we do not yet know the multiplicity of each distinct prime factor.

To go all the way to the natural numbers, we must account for the missing multiplicity. This is handled
by the second part of the constructive prime generator, contained in the following renormalisation procedure: Consider a prime ruler with unit $p_k$ and assign it the role as a unit ruler. Then construct a new ruler with unit $p_k$, measured in terms of the assigned unit ruler. The resulting ruler in turn takes on the role as unit ruler and the process is repeated indefinitely. For any $k \geq 1$, the outcome is a set of  \emph{renormalised prime rulers} with units $p_k^2$, $p_k^3$, $p_k^4$, $\dots$, measured in terms of the original unit ruler with unit 1. The case of $k=1$ is illustrated in Figure~1B. Similarly, from the numerical perspective, we generate for each $k$ the \emph{renormalised prime sequences} with periods ${p_k^2, p_k^3, p_k^4, \dots}$ and so on. Note, however, that the repeating elements in each of these sequences are still 1 and $p_k$, as shown in \tabref{tab:table1}B for $k=1$.
%
%
%
%
%

Combining the two parts of the constructive prime generator, we are now able to measure any distance on the original unit ruler, beyond one unit, in terms of prime rulers and renormalised prime rulers, or, more to our interest, define any natural number greater than 1 in terms of the prime sequences and their renormalised counterparts. Essentially, we have arrived at a complete bottom-up picture of the multiplicative architecture underlying the natural numbers, an architecture entirely composed of periodic sequences, as depicted in  \tabref{tab:table2}.
\begin{table}
\caption{\small \label{tab:table2} The multiplicative architecture underlying the natural numbers. The set of natural numbers can be viewed (or defined) as the multiplicative combination of the set of prime periodic sequences and their renormalised counterparts.}
\centering
\resizebox{0.85\columnwidth}{!}{%
\begin{tabular}{ccccccccccccccccccc}
\obsgray{$0$} & 1 & 2 & 3 & 4 & 5 & 6 & 7 & 8 & 9 & 10 & 11 & 12 & 13 & $\cdots$ & $p_k$ & $\cdots$& $p_k^2$ & $\cdots$\vspace{1mm}\\
\hline  \\ [-1mm]
\obsgray{$1$} & 1 & $1$ & 1 & $1$ & 1 & $1$ & 1 & $1$ & 1 & $1$ & 1 & $1$ & 1 & $\cdots$ & 1 &  $\cdots$& 1 &  $\cdots$\\
\obsgray{$p_1$} & 1 & $p_1$ & 1 & $p_1$ & 1 & $p_1$ & 1 & $p_1$ & 1 & $p_1$ & 1 & $p_1$ & 1 & $\cdots$ & 1 &  $\cdots$& 1 &  $\cdots$\\
\obsgray{$p_1$} & 1 & $1$ & 1 & $p_1$ & 1 & $1$ & 1 & $p_1$ & 1 & $1$ & 1 & $p_1$ & 1 & $\cdots$ & 1 &  $\cdots$& 1 &  $\cdots$\\
\obsgray{$p_1$} & 1 & $1$ & 1 & 1 & 1 & $1$ & 1 & $p_1$ & 1 & $1$ & 1 & 1 & 1 & $\cdots$ & 1 &  $\cdots$& 1 &  $\cdots$\\
\vdots & \vdots & \vdots & \vdots & \vdots & \vdots & \vdots & \vdots & \vdots & \vdots & \vdots & \vdots & \vdots & \vdots & $\cdots$ 	& \vdots & $\cdots$& \vdots &  $\cdots$\\
\obsgray{$p_2$} & 1 & 1 & $p_2$ & 1 & 1 & $p_2$ & 1 & 1 & $p_2$ & 1 & 1 & $p_2$ & 1 & $\cdots$ & $1$ & $\cdots$ & 1 &  $\cdots$\\
\obsgray{$p_2$} & 1 & 1 & $1$ & 1 & 1 & $1$ & 1 & 1 & $p_2$ & 1 & 1 & $1$ & 1 & $\cdots$ & $1$ & $\cdots$ & 1 &  $\cdots$\\
\obsgray{$p_2$} & 1 & 1 & $1$ & 1 & 1 & $1$ & 1 & 1 & 1 & 1 & 1 & $1$ & 1 & $\cdots$ & $1$ & $\cdots$ & 1 &  $\cdots$\\
\vdots & \vdots & \vdots & \vdots & \vdots & \vdots & \vdots & \vdots & \vdots & \vdots & \vdots & \vdots & \vdots & \vdots & $\cdots$ 	& \vdots & $\cdots$& \vdots &  $\cdots$\\
\obsgray{$p_k$} & 1 & 1 & 1 & 1 & 1 & 1 & 1 & 1 & 1 & 1 & 1 & 1 & 1 & $\cdots$ & $p_k$ &$\cdots$& $p_k$ &  $\cdots$\\ 
\obsgray{$p_k$} & 1 & 1 & 1 & 1 & 1 & 1 & 1 & 1 & 1 & 1 & 1 & 1 & 1 & $\cdots$ & 1 &$\cdots$& $p_k$ &  $\cdots$\\ 
\obsgray{$p_k$} & 1 & 1 & 1 & 1 & 1 & 1 & 1 & 1 & 1 & 1 & 1 & 1 & 1 & $\cdots$ & 1 &$\cdots$& 1 &  $\cdots$\\ 
\vdots & \vdots & \vdots & \vdots & \vdots & \vdots & \vdots & \vdots & \vdots & \vdots & \vdots & \vdots & \vdots & \vdots & $\cdots$ 	& \vdots & $\cdots$ & \vdots & $\cdots$
\end{tabular}%
}
\end{table}

As stated earlier, the prime numbers are frequently presented as the multiplicative building blocks of the natural numbers. But our exposure of the architecture underneath the natural numbers suggests the alternative and in some sense more precise statement: \emph{The multiplicative building blocks of the natural numbers are the prime sequences and their corresponding renormalised sequences.} This statement is perhaps a trivial reformulation, but indirectly it emphasises the fact that understanding the properties of the prime numbers in relation to the natural numbers is fundamentally about understanding the properties of an ever increasing combination of periodic sequences. As we detail in the next section, this insight can be exploited to construct a random model of the primes built on the same structure as the primes. 

With the constructive prime generator completed, it remains to properly define the natural numbers in terms of our sets of periodic sequences. We start by expressing all prime sequences and renormalised prime sequences on the form $\rho_{i,j}(n)$, which refers to the sequence with period equal to $p_i^j$.
In addition we denote the unit sequence by $\rho_0(n)$. The natural number $n$ can therefore be stated either as the infinite product
\begin{align}
	n
	:= 
	\rho_0(n) \cdot \prod_{i=1}^\infty \prod_{j=1}^\infty\rho_{i,j}(n) 
	= 
	\prod_{i=1}^\infty \prod_{j=1}^\infty\rho_{i,j}(n),
\label{eq:canonical1}
\end{align}
or, as the finite product 
\begin{align}
	n
	:= 
	\prod_{i\in I_n} \prod_{j=1}^{l_i(n)}\rho_{i,j}(n), 
\label{eq:canonical2}
\end{align}
where $I_n$ denotes the set of indices of the distinct prime factors of $n$ and $l_i(n)$ counts the multiplicity of $p_i$ in $n$. Note that we are deliberately being pragmatic with our notation here, as $n$ appears on both sides of these equations (as well as our use of integer indices). One could in fact think of $n$ when it appears on the right hand side as just a name for, or pointer to, the respective position in the unit sequence. On the left hand side, $n$ is assigned the meaning of a product of primes. As a curiosity, this view resembles what Goethe wrote already centuries ago~\citep{Newman:1956}: \emph{Two times two is not four, but it is just two times two, and that is what we call four for short.} 

Obviously, \eqref{eq:canonical1} and \eqref{eq:canonical2} are equivalent to the usual canonical expressions 
$${n =  \prod_{i=1}^{\infty} p_i^{m_i}} \quad \textrm{and} \quad {n =  \prod_{i\in I_n} p_i^{\l_i(n)}},$$
respectively, where at most a finite number of ${m_i}$ are positive integers, and the remaining are zero. The possible advantage of \eqref{eq:canonical1} and \eqref{eq:canonical2}, if only pedagogically, is that these expressions explicitly emphasise that the natural numbers are built up multiplicatively from periodic sequences. In particular, we have from \eqref{eq:canonical1} that the sequence of natural numbers is the limit of the periodic sequence ${\rho_0(n) \cdot \prod_{i=1}^m \prod_{j=1}^k\rho_{i,j}(n)}$ as $m, k$ go towards infinity. In other words, we obtain the natural numbers by iterating the constructive prime generator indefinitely.


The uniqueness of $n$ is of course guaranteed by the fundamental theorem of arithmetic, which states:
\begin{theorem}[Fundamental theorem of arithmetic] 
Every natural number greater than 1 can be expressed in exactly one way as a product of primes, apart from rearrangement of factors.
\end{theorem}
Standard proofs of this theorem, see for example \citep[p. 11]{Tenenbaum:2015}, assumes the number theoretical definition of primes, where the primes are defined in terms of the natural numbers. Our starting point, however, is reversed; we first generate the primes and subsequently define the natural numbers. This leads to an alternative proof of the fundamental theorem of arithmetic: 
\vspace{2mm}
\begin{proof} By our definition of the natural numbers, \eqref{eq:canonical1} or \eqref{eq:canonical2}, any natural number greater than 1 can be expressed as a product of one or more primes. Uniqueness follows directly from the fact that the natural numbers are constructed out of periodic sequences: Assume that the position $n$ in the unit sequence is the first position measured by the product ${\prod_{i\in I} \prod_{j=1}^{l_i}\rho_{i,j}(n)}$, where $I$ is a finite subset of the positive integers and $1 \leq l_i < \infty$. This product is periodic with respect to the position $n$, with period equal to the product itself. The only other positions also measurable by this product are therefore of the form $k n$, where $k>1$. But, again by definition, every $k$ is a product of at least one prime, so $n$ is the only position fully measured by ${\prod_{i\in I} \prod_{j=1}^{l_i}\rho_{i,j}(n)}$.
\end{proof}
\vspace{2mm}

Up to this point we have directed our attention as to how (and why) the prime numbers should be conceived of as periodic structures rather than as single numbers. Let us next see how this perspective is useful for understanding the structure of the primes in the natural numbers and how the apparent random nature of the primes is a direct consequence of the prime sequences themselves being recursively generated from preceding prime sequences.

\section{The structure and randomness of the prime numbers}
\label{sec:structure}

In the typical presentation of the sieve of Eratosthenes, as for example in \citep[p. 1]{Friedlander:2010}, one usually considers a given value $x$, and then proceeds to find the primes smaller or equal to $x$ by removing composite numbers, starting with those that are composites of $p_1=2$, then $p_2=3$, and so on. Since any composite which has prime factors all larger than $\sqrt{x}$ must necessarily exceed $x$, all natural numbers remaining after sieving by the prime numbers smaller or equal to $\sqrt{x}$ are primes, and the sieving process is therefore completed at this point. For example, with $x=45$ we only need to remove composites of $2, 3,$ and $5$ to locate all primes below $x$, as illustrated in \tabref{tab:table3}.
\begin{table}[t]
\caption{\small \label{tab:table3}Illustration of the sieve of Eratosthenes. To find all primes smaller or equal to  $x$, we start by removing composites of the first prime $p_1=2$, and continue removing composites of all consecutive primes smaller or equal to $\sqrt{x}$. In this example, all primes below $x=45$ are found by removing composites of 2, 3, and 5. 
}
\centering
\resizebox{0.8\columnwidth}{!}{%
\begin{tabular}{ccccccccccccccc}
1 & 2 & 3 & \sieved{4} & 5 & \sieved{6} & 7 & \sieved{8} & \sieved{9} & \sieved{10} & 11 & \sieved{12} & 13 & \sieved{14} & \sieved{15}\\
\sieved{16} & 17 & \sieved{18} & 19 & \sieved{20} & \sieved{21} & \sieved{22} & 23 & \sieved{24} & \sieved{25} & \sieved{26} & \sieved{27} & \sieved{28} & 29 & \sieved{30}\\
31 & \sieved{32} & \sieved{33} & \sieved{34} & \sieved{35} &  \sieved{36} & 37 & \sieved{38} & \sieved{39} & \sieved{40} & 41 & \sieved{42} & 43 & \sieved{44} & \sieved{45} 
\end{tabular}%
}
\end{table}

This way of presenting the sieve of Eratostenes, however, leaves out essential information about the structure of the primes. 
To see why, let us initially remove the focus on the given value $x$, and assume we sieve to infinity in each step. Then we can make the observation that the smallest composite removed in the $k$th sieve step is $p_{k}^2$, while the smallest composite remaining after the $k$th sieve step is $p_{k+1}^2$. In other words, in the $k$th sieve step we complete the sieving of the integer set ${{s_k:=\{p_k^2, \dots, p_{k+1}^2-1\}}}$, which
holds true even for $k=0$ if we define $p_0$ as 1. 

Now, let us interpret this observation in terms of the constructive prime generator of the previous section. Only the first part of the generator is relevant, as this provides us with the primes, so for convenience, let us write $\rho_i(n):= \rho_{i,1}(n)$. Then each integer set $s_k$ has the specific property that all its primes are recursively generated from the unit sequence $\rho_0(n)$ and the previous prime sequences $\rho_i(n)$, $1\leq i \leq k$. An example is shown for $s_3$ in \tabref{tab:table4}, where the generated primes appear in the columns with all 1s.
\begin{table}[t]
\caption{\small \label{tab:table4} The recursive structure of the primes in the integer set ${s_3:=\{p_3^2, \dots, p_{4}^2 -1\}}$. Primes in $s_3$ are generated recursively from the unit sequence $\rho_0(n)$ and the $3$ first prime sequences $\rho_1(n)$, $\rho_2(n)$, and $\rho_3(n)$ for all columns containing only 1s.}
\centering
\resizebox{0.8\columnwidth}{!}{%
\begin{tabular}{c|ccccccccccccccc|}
\centering
\obsgray{$n$} & 25 & 26 & 27 & 28 & 29 & 30 & 31 & 32 & 33 & 34 & 35 & 36 & 37 & $\cdots$ & 48\\
\hline  \\ [-3mm]
\obsgray{$\rho_0(n)$}  & 1 & 1 & 1 & 1 & \bf{1} & 1 & \bf{1} & 1 & 1 & 1 & 1 & 1 & \bf{1} &$\cdots$ & 1\\
\obsgray{$\rho_1(n)$}  & 1 & $p_1$ & 1 & $p_1$ & \bf{1} & $p_1$ & \bf{1} & $p_1$ & 1 & $p_1$ & 1 & $p_1$ & \bf{1} &$\cdots$ & $p_1$\\
\obsgray{$\rho_2(n)$}  &1 & 1 & $p_2$ & 1 & \bf{1} & $p_2$ & \bf{1} & 1 & $p_2$ & 1 & 1 & $p_2$ & \bf{1} &$\cdots$ & $p_2$\\
\obsgray{$\rho_3(n)$}  &$p_3$ & 1 & 1 & 1 & \bf{1} & $p_3$ & \bf{1} & 1 & 1 & 1 & $p_3$ & 1 & \bf{1} &$\cdots$ & 1
\end{tabular}%
}
\end{table}

Combining this property---that the primes in any integer set $s_k$ are generated recursively from the $k$ first prime sequences---with the fact that the integer sets $s_k$ make up a complete subdivision of the natural numbers, we arrive at a clear view of how the recursive structure of the primes is manifested in the natural numbers (Figure~2). Building on this perspective, we will spend the reminder of this section deepening our understanding of the distribution of primes in the natural numbers.  

\begin{figure}[t]
\centering
\includegraphics[width=0.8\textwidth]{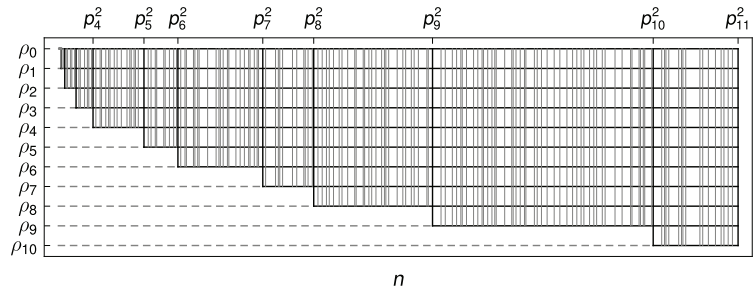}
\caption{\small 
The recursive structure of the primes from the perspective of the constructive prime generator. 
The dashed horizontal lines represents the prime sequences $\rho_i(n)$, including the unit sequence $\rho_0(n)$. The vertical gray lines show the positions where the product $\prod_{0\leq i \leq k} \rho_i (n) =1$ for $n\in s_k$, and hence, where new primes must be generated. The view promoted by this diagram is that the distribution of primes in each integer set $s_k$ is recursively generated from the unit sequence and the preceding prime sequences $\rho_i(n)$, $0\leq i \leq k$.}
\label{fig:structure}
\end{figure}

\subsubsection{The characteristic function of prime numbers}
A straightforward way to couple the distribution of primes to their structure in the natural numbers, is by expressing the characteristic function of primes, $\mathbf{1}_\mathcal{P}(n)$, in terms of the prime sequences $\rho_k(n)$, $k\geq 0$ (unit sequence included). Rather than $\rho_k(n)$ though, it is convenient to define the equivalent function $\kappa_k(n)$, which is identical to $\rho_k(n)$ whenever this takes the value one and zero otherwise. 
We now apply this definition to write the characteristic function of primes as 
\begin{align}
{\mathbf{1}_\mathcal{P}(n) = \prod_{k=0}^{\pi(\sqrt{n}) } \kappa_k(n)},
\label{eq:characteristicfunction}
\end{align}
where $\pi(x) :=  \sum_{i=1}^{[x]} \mathbf{1}_\mathcal{P}(i) - 1$ is the prime counting function. Note that ${\mathbf{1}_\mathcal{P}(1) = 1}$ by construction, which explains the subtracting factor of $-1$ in this definition of $\pi(x)$. For a formulation in terms of standard number theoretical notation, see Supplementary information (SI) Standard notation.
%
%

This formulation of the characteristic function defines a recurrence relation for the distribution of primes in the natural numbers, and reinforces the view that the primes in the $k$th integer set $s_k$ should be viewed as generated recursively from the $k$ first prime sequences. It is instructive to observe that the right hand side of \eqref{eq:characteristicfunction} contains a multiplicative combination of deterministic periodic sequences, while the left hand side provides the distribution of the individual primes, which we know behaves randomlike. It is therefore clear that the apparent randomness of the primes must be entirely explainable in terms of an ever growing combination of periodic sequences. Working directly on the primes can be complicated, however, but as we shall see next, our formulation of the characteristic function $\mathbf{1}_\mathcal{P}(n)$ immediately suggests a random model that can be expressed and analysed theoretically and experimentally, with implications also for our understanding of the distribution of primes.

\subsubsection{A random model of the primes}

The idea behind our random model is simply to shift the structure of the primes randomly relative to the natural numbers. This is easily visualised by imagining the structure shown in Figure~2 placed at random positions in the natural numbers. The "primes" in this model then appear in the columns containing only 1s, as is the case for the primes themselves. Theoretically, we define the random model (RM) in terms of the characteristic function 
\begin{align}
	\mathbf{1}_{\mathcal{P}_{\operatorname{RM}}}(n)
	:= 
	\prod_{k=0}^{\pi(\sqrt{n})} \kappa_k(n+a),
\label{eq:RM}	
\end{align}
where $a$ is a random integer and ${\mathcal{P}_{\operatorname{RM}}}$ is the set of "primes" generated by a realisation of the RM. Technically, $a$ is defined so that $a \equiv b_k \pmod{  p_k}$ with probability $1/p_k$ for all ${b_k \in \{0, 1, \dots, p_k-1\}}$ and ${0\leq k \leq \pi(\sqrt{n})}$. Notice now the strong resemblance between $\mathbf{1}_{\mathcal{P}_{\operatorname{RM}}}(n)$ and $\mathbf{1}_\mathcal{P}(n)$: For a given realisation of $a$, the characteristic function $\mathbf{1}_{\mathcal{P}_{\operatorname{RM}}}(n)$ is completely deterministic, as is $\mathbf{1}_\mathcal{P}(n)$, and they both derive from the same underlying structure of periodic sequences. As a result, the RM only allows for "prime patterns" that occur in the real primes, such as twin primes, or more generally, prime k-tuples. In contrast, models like Cramer's random model \citep{Cramer:1936ud} or Hawkins' random sieve \citep{Hawkins:1957} have no restrictions on what patterns can occur, though modified models exist that enforce restrictions on obtainable patterns in the long term limit, as for example \citep[p. 66]{Tenenbaum:2015}. To aid in the understanding of how the RM behaves compared to Cramer's model, we present a visual comparison of both models in SI Figure S1. 

%
%
%
%

%
%

To examine the relation between the primes and the RM, we consider the sum across ${\mathbf{1}_{\mathcal{P}_{\operatorname{RM}}}(n)}$ in terms of the counting function 
${
	\pi_{\operatorname{RM}}(x) 
	:=  
	\sum_{i=1}^{[x]} \mathbf{1}_{\mathcal{P}_{\operatorname{RM}}}(i) - 1
}$.
By introducing ${W(x):= \prod_{p\leq x} \left(1-1/p\right)}$, the expected value of $\pi_{\operatorname{RM}}(x)$ is simply
\begin{align*}
	\mathbf{E}\left[ \pi_{\operatorname{RM}}(x) \right]
	=  
	\sum_{n=2}^{[x]} W(\sqrt{n})
	\sim 
	2 \operatorname{e}^{-\gamma} \operatorname{li}(x),
\end{align*}
where the asymptotic equality follows from Merten's product theorem~\citep[p. 19]{Tenenbaum:2015}. Here ${\operatorname{li}(x)}$ is the logarithmic integral $\int_2^x dt/\log t$, while $\gamma$ is the Euler-Mascheroni constant. The prime number theorem on the other hand states that $\pi(x) \sim \operatorname{li}(x)$, so it follows that, as $n \to \infty$, the two counting functions $\pi_{\operatorname{RM}}(x)$ and $\pi(x)$ satisfy the relation
\begin{align*}
	\mathbf{E}\left[ \pi_{\operatorname{RM}}(x) \right]
	\sim 
	2 \e^{-\gamma} \pi(x).
\end{align*}

Since the set of primes corresponds to only one out of all possible realisations of the RM, it is nothing strange about the fact that the actual asymptotic of $\pi(x)$ deviates from the expected value of $\pi_{\operatorname{RM}}(x)$. Indeed,  the reason why these two means differ follows from the simple fact that the product of the unique prime factors of $n$ can never exceed $n$, a constraint we did not account for in the RM. In other words, when we earlier replaced  $\rho_k(n)$ by $\kappa_k(n)$ in order to define $\mathbf{1}_\mathcal{P}(n)$ and $\mathbf{1}_{\mathcal{P}_{\operatorname{RM}}}(n)$, we did not incorporate the fact that $\prod_{k=0}^{\pi(\sqrt{n})} \rho_k(n) \leq n$. If we impose this constraint on the RM, considering only the subset of realisations that satisfy ${\prod_{k=0}^{\pi(\sqrt{n})} \rho_k(n+a) \leq  n}$  for all values of $n$, we observe that the two means appear to be the same asymptotically (Figure~3), lending support to the conjecture: 
\begin{conjecture} 
Let ${\operatorname{RM}}_c$ denote the constrained ${\operatorname{RM}}$ as explained above. The expected value of $\pi_{\operatorname{RM}_c}(x)$ satisfies
\begin{align*}	
	{\mathbf{E}\left[ \pi_{\operatorname{RM}_c}(x) \right]
	\sim 
	{\operatorname{li}}(x)}.
\end{align*}
\end{conjecture}

To strengthen this conjecture, note that for a fixed value of $n$, and assuming only values of $a$ so that ${\prod_{k=0}^{\pi(\sqrt{n})} \rho_k(n+a) \leq  n}$, the expected value of $\prod_{k=0}^{\pi(\sqrt{n})} \kappa_k(n+a)$ can be written as
\begin{align*}	
	\mathbf{E}\left [\prod_{k=0}^{\pi(\sqrt{n})} \kappa_k(n+a)\right] 
	= 
	\sum_{ \substack{ d \mid P(\sqrt{n}) \\ d \leq n  }}^{} \frac{\mu(d)}{d},
\end{align*}
where $\mu (d)$ is the M\"obius function and $P(x)$ denotes the product of all primes smaller or equal to $x$. For the latter expression one can then prove the following result (proof provided in SI Proposition S3.1):
\begin{proposition}
The sum $\sum_{ \substack{ d \mid P(\sqrt{n}) \\ d \leq n  }}^{} \frac{\mu(d)}{d}$ satisfies the asymptotic relation
\begin{align*} 	
	\sum_{ \substack{ d \mid P(\sqrt{n}) \\ d \leq n  }}^{} \frac{\mu(d)}{d}
	\sim 
	\frac{1}{\log n}.
\end{align*}
\end{proposition}

Now, if we for each value of $n \geq 2$ draw a new value of $a$, say $a_n$, each satisfying  
${\prod_{k=0}^{\pi(\sqrt{n})} \rho_k(n+a_n) \leq  n}$,
the result is a set of uncorrelated random variables whose expected sum obeys the same asymptotic as the prime number theorem,
\begin{align*}
	\mathbf{E}\left [ \sum_{n \leq x} \prod_{k=0}^{\pi(\sqrt{n})} \kappa_k(n+a_n) \right]
	\sim
	\sum_{n \leq x} \frac{1}{\log n}
	\sim 
	\li(x).
\end{align*}
Necessarily, any realisation of the $\operatorname{RM}_c$ in terms of $\pi_{\operatorname{RM}_c}(x)$ lies in a subspace of the outcomes possible by  $\sum_{n \leq x} \prod_{k=0}^{\pi(\sqrt{n})} \kappa_k(n+a_n)$, suggesting $\pi_{\operatorname{RM}_c}(x)$ shares the same expected asymptotic mean. Thus, not only does the $\operatorname{RM}$ completely account for the structure of the primes, it presumably also accounts for the correct asymptotic density of the primes if we consider the $\operatorname{RM}_c$. 

\begin{example}
The sample space of ${\pi_{\operatorname{RM}_c}(x)}$ both grow and shrink in size with increasing $x$ (SI Table S1), in contrast to the full sample space of ${\pi_{\operatorname{RM}}(x)}$, which branches out by a factor of $p_k$ each time $x$ enters an interval $[p_k^2,p_{k+1}^2)$. For a numerical example, consider $x=p_{41}^2-1$. The sample space of ${\pi_{\operatorname{RM}}(x)}$ then has size ${\prod_{i \leq 40} p_i > 1.6 \times 10^{68}}$, while ${\pi_{\operatorname{RM}_c}(x)}$ contains only 88 elements in all. As conjectured above, the realisations of ${\pi_{\operatorname{RM}_c}(x)}$ all lie close to the asymptotic mean $\operatorname{li}(x)$, while a corresponding set of random realisations of ${\pi_{\operatorname{RM}}(x)}$ cluster around the expected value ${\mathbf{E}\left[ \pi_{\operatorname{RM}}(x) \right]}$ (Figure~3A). Upon a closer look, all realisations of ${\pi_{\operatorname{RM}_c}(x)}$ are strongly correlated, so even the mean value fluctuates accordingly, and there are no outliers in this constrained sample space---at least on this scale of numerical evidence  (Figure~3B). 
This suggests the even stronger conjecture that 
$
	{\pi_{\operatorname{RM}_c}(x)
	\sim 
	{\operatorname{li}}(x)}
$
for all realisations of the $\operatorname{RM}_c$.
\end{example}

Having explained the average behaviour of the RM and its relation to the primes and the ${\operatorname{RM}_c}$, let us now turn our attention towards the fluctuations around the average behaviour and an explanation of the randomlike nature of the primes. 
%

\begin{figure}[t]
\centering
\includegraphics[width=\textwidth]{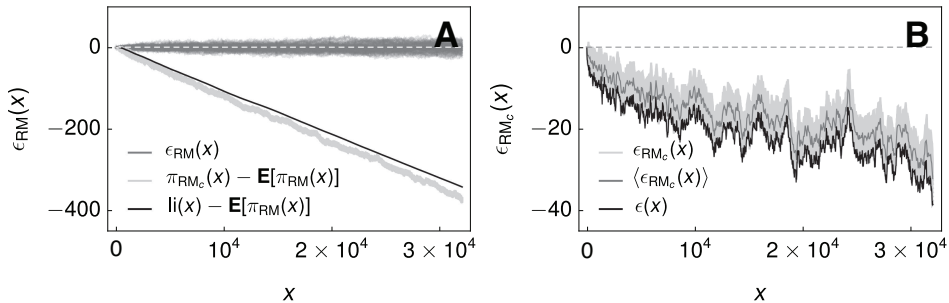}
\caption{\small
Realisations of the $\operatorname{RM}$ for $k=40$ sieve steps. A: The dark gray curves show 
${\epsilon_{\operatorname{RM}}(x):={\pi_{\operatorname{RM}}(x)} - \mathbf{E}[{\pi_{\operatorname{RM}}(x)}]}$ 
for 88 (out of more than ${{1.6 \times 10^{68}}}$) random realisations of the $\operatorname{RM}$. Likewise, the light gray curves show all 88 realisations of the ${\operatorname{RM}_c}$ present at $x=p_{41}^2-1$, all lying close to 
${\operatorname{li}(x) - \mathbf{E}[{\pi_{\operatorname{RM}}(x)}]}$
(black). B: The same realisations of the ${\operatorname{RM}_c}$, now plotted in terms of 
${\epsilon_{\operatorname{RM}_c}(x):={\pi_{\operatorname{RM}_c}(x)} - \operatorname{li}(x)}$.
Also shown is the mean value ${\langle \epsilon_{\operatorname{RM}_c}(x)\rangle}$ (darker gray), as well as $\epsilon(x):=\pi(x)-\operatorname{li}(x)$ (black). 
All realisations of ${\epsilon_{\operatorname{RM}_c}(x)}$ are strongly correlated, suggesting $\pi_{\operatorname{RM}_c}(x) \sim \operatorname{li}(x)$ to be true for all realisations of the ${\operatorname{RM}_c}$.
}
\label{fig5}
\end{figure}

\subsection{The randomness of primes}

The apparent randomness of the prime numbers manifests itself most clearly in central limit behaviour, one example of which is the famous central limit theorem by Erd\H{o}s-Kac~\citep{Erdos:1940}, stating that the number of distinct prime factors of $n$ is normal distributed with mean and variance asymptotically equal to $\log\log n$. As we shall see here, central limit behaviour in the primes---and thus randomlike behaviour---is a natural consequence of how the distribution of primes is recursively built up from previous prime periodic sequences. 
%
%
%
%
%
%
%
%
%

To begin with, let us make the observation that the prime periodic sequences $\kappa_k(n)$ and any finite multiplicative combination of these belong to a category of sequences that are known as subrandom (or low-discrepancy) sequences. Such sequences can stem from correlated random processes or even be completely deterministic, but, in general, they have a correlation structure that forces the underlying distribution to be sampled more efficiently as compared to an uncorrelated random process with the same expected value. In particular, a subrandom sequence is dominated by negative correlations between its elements, and a sum over such a sequence will therefore have lower variance as compared to a corresponding sum deriving from a fully random process.
%
%
%
%
%
%

For a concrete and relevant example, that also brings in the aspect of central limit behaviour, consider the multiplicative combination of prime periodic sequences given by 
${K_k(n):=\prod_{i=0}^k \kappa_i(n)}$ and the corresponding random sum 
$${S_k(h) :=\sum_{n=1}^{h} K_k(n+ a)},$$ 
where $a$ is a random integer as explained in connection with $\eqref{eq:RM}$. The sequence $K_k(n)$ is itself periodic, with period equal to $\prod_{i=1}^k p_k$, and regularity is trivially maintained on this scale. On a local scale, however, $K_k(n)$ becomes increasingly irregular as $k$ gets larger, extending the range of possible outcomes of ${S_k(h)}$. The result is that the distribution of ${S_k(h)}$ eventually approaches a normal distribution, for 
an appropriate range of $h$ that depends on $k$.
This fact was proven in \citep{Montgomery:1986} and later improved upon in \citep{Montgomery:2004de}. 
%
The variance of $S_k(h)$ was first derived in \citep{hausman:1973} and it is a simple matter to prove that for $h>1$ the variance is strictly less than $h W(p_k)(1-W(p_k))$, corresponding to the variance when the sum stems from a random process. The proof is similar to that of \propositionref{thm:covar}, which is provided in SI Proposition S4.1. 

\begin{remark}
The lesson to take away from this example is that we attain apparent randomness---in terms of central limit behaviour---simply by multiplicatively combining a finite set of prime periodic sequences. The exact same mechanism underlies the randomness of the primes, a consequence of the distribution of primes in each integer set $s_k$ being recursively generated by the multiplicative combination of $k$ prime periodic sequences.  
\end{remark}  
%
%
%
%
%
%
%
%
%
%
%
%

With this background, let us now turn to the RM. 
For later comparison with the primes it is useful to write the variance of 
${\pi_{\operatorname{RM}}(x)}$
on the form 
\begin{align}
\label{eq:RMvar}
	\operatorname{Var}(\pi_{\operatorname{RM}}(x))
	= 
	\sum_{i=1}^{[x]} \operatorname{Var}(\mathbf{1}_{ \mathcal{P}_{\operatorname{RM}} }(i))
	+
	2 \sum_{i=1}^{[x]} \sum_{i<j}^{[x]} 
	\operatorname{Cov}
	\left(
	\mathbf{1}_{ \mathcal{P}_{\operatorname{RM}} }(i),
	\mathbf{1}_{ \mathcal{P}_{\operatorname{RM}} }(j)
	\right).\nonumber
\end{align}
The first sum on the right hand side---the sum of variances---equals
$${
\sum_{i=1}^{[x]} W(\sqrt{i})(1-W(\sqrt{i})),
}$$
which is the variance we would obtain for a sum over uncorrelated random variables with the same expected values. 
The second sum---the sum of covariances---can be derived analytically by
%
%
%
%
%
%
%
generalising the derivation in \citep{hausman:1973} (SI Covariance Expression), from which we obtain the following result: 
%
%
%
\begin{proposition} 
The covariance  
$ 
	\operatorname{Cov}
	\left(
	\mathbf{1}_{ \mathcal{P}_{\operatorname{RM}} }(i),
	\mathbf{1}_{ \mathcal{P}_{\operatorname{RM}} }(j)
	\right)
$
satisfies the relation
\begin{align*}
	\sum_{i<j}^{[x]} 
	&\operatorname{Cov}
	\left(
	\mathbf{1}_{ \mathcal{P}_{\operatorname{RM}} }(i),
	\mathbf{1}_{ \mathcal{P}_{\operatorname{RM}} }(j)
	\right)
	\leq 0, 
\end{align*}
where the inequality is strict for $i \geq 25$. 
\label{thm:covar}
\end{proposition}
While this proposition (see proof in SI Proposition S4.1) guarantees that the variance of $\pi_{\operatorname{RM}}(x)$
is strictly smaller than that deriving from a random process, the numerical evidence in Figure 4A in fact suggests the stronger conjecture 
\begin{align*}
	\limsup_{x \to \infty}
	\sum_{i=1}^{[x]} \sum_{i<j}^{[x]} 
	\operatorname{Cov}
	\left(
	\mathbf{1}_{ \mathcal{P}_{\operatorname{RM}} }(i),
	\mathbf{1}_{ \mathcal{P}_{\operatorname{RM}} }(j)
	\right)	= 
	-\infty.
\end{align*}
%

\begin{figure}[b]
\centering
\includegraphics[width=\textwidth]{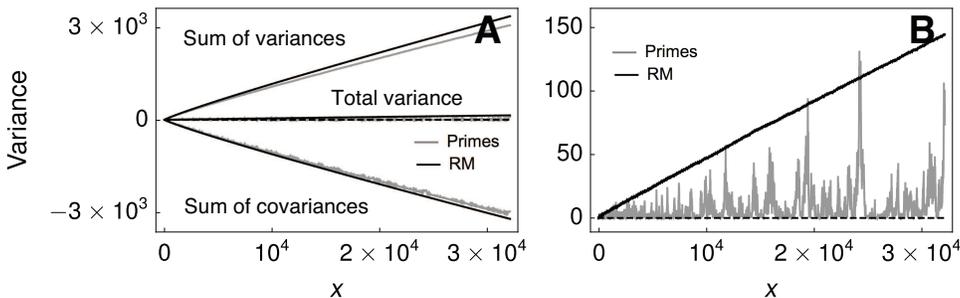}
\caption{\small
Comparing $\operatorname{Var}[\pi_{\operatorname{RM}}(x)]$ and $|\pi(x)-\operatorname{Ri}(x)|^2$.
A: In black, 
the variance of the $\operatorname{RM}$ (left side of \eqref{eq:RMvar}) plotted separately and in terms of its split components (right side of \eqref{eq:RMvar}). Likewise, in gray, the squared error for the primes and its split components \eqref{eq:primevar}. Qualitatively, the RM provides an accurate description of the fluctuating behaviour of the primes. B: Close-up view, showing only the variance of the $\operatorname{RM}$ and the squared error for the primes. 
}
\end{figure}


%

%
%
%

What now for the $\operatorname{RM}_c$ and the primes? 
In the absence of a rigorous proof, we take a numerical approach to investigating the variance of $\pi_{\operatorname{RM}_c}(x)$ with respect to the asymptotic mean $\operatorname{li}(x)$. Given the observed strong correlations between the different samples of the $\operatorname{RM}_c$ (Figure 3B), it should suffice to consider one realisation, which we choose to be $\pi(x)$.
%
%
%
%
%
%
Numerically, the cleanest result is obtained by employing the Riemann function
\begin{align}
\operatorname{Ri}(x):=\sum_{m=1}^{\infty} \frac{\mu(m)}{m}\operatorname{li}(x^{1/m}) 
\end{align}
 as the expected value for $\pi(x)$, as ${\operatorname{Ri}}(x)$ on average is a better estimate than $\operatorname{li}(x)$ for $\pi(x)$ \citep[p. 105-106]{Ingham:1932}. These two estimates are asymptotically equal, however, and the result holds also for $\operatorname{li}(x)$ when $x$ is large (SI Figure S2). 
 
 Specifically, we want to compare the terms in \eqref{eq:RMvar} 
with the corresponding terms in 
\begin{align}
\label{eq:primevar}
	|\pi(x) - \operatorname{Ri}(x)|^2
	= 
	\sum_{i=1}^{[x]} \epsilon_i^2
	+
	2 \sum_{i=1}^{[x]} \sum_{i<j}^{[x]} \epsilon_i \epsilon_j, 
\end{align}
where  
${\epsilon_i:= \mathbf{1}_{\mathcal{P}}(i) -(\operatorname{Ri}(i) - \operatorname{Ri}(i-1))}$. 
Firstly, we observe that the two terms on the right hand side of \eqref{eq:primevar}---the sum of variances and the sum of covariances---behave similarly to the corresponding terms in \eqref{eq:RMvar} (Figure 4A).
For the sum of variances, it is easy to show that the ratio of the term in \eqref{eq:RMvar} to the term in \eqref{eq:primevar} will tend to  $2\e^{-\gamma}$. While the exact relation between the two sums of covariances awaits a theoretical investigation, we find that numerically, the RM provides a qualitatively accurate description of the fluctuating behaviour of the primes, and as for the RM, we conjecture that 
\begin{align*}
	\limsup_{x \to \infty}
	\sum_{i=1}^{[x]} \sum_{i<j}^{[x]} \epsilon_i \epsilon_j 
	= 
	-\infty.
\end{align*}
Secondly, comparing the left hand terms in \eqref{eq:RMvar} and \eqref{eq:primevar}, we observe that the squared error in the case of the primes---and presumably the variance of the $\operatorname{RM}_c$---on average is smaller than the variance of the RM  (Figure 4B). 
This is expected, and can be anticipated from the heuristic argument that, effectively, the $\operatorname{RM}_c$ amounts to picking samples biased towards one side of a normal distribution, thereby reducing the variance as compared to the original distribution.  

The steadily increasing variance of the RM (Figure 4B) reveals that $\pi_{\operatorname{RM}}(x)$ exhibits larger fluctuations on average as $x$ increases. One way of viewing this is that the correlations between sequence elements diminish as $x$ grows. More precisely, one can prove that (see SI Proposition S4.2)
\begin{proposition}
Let $\operatorname{Corr}[x,y]$ be the correlation function between two variables $x$ and $y$. Then we have that the RM in terms of the indicator function $\mathbf{1}_{\mathcal{P}_{\operatorname{RM}}}(m)$ satisfies the relation
\begin{align*}
\lim_{m,n\to \infty} \operatorname{Corr}[\mathbf{1}_{\mathcal{P}_{\operatorname{RM}}}(m), \mathbf{1}_{\mathcal{P}_{\operatorname{RM}}}(n)] =0,
\end{align*}
for any fixed distance $n-m>0$.
\end{proposition}
%
The growth in fluctuations is exactly what happens in the primes as well. In fact, we have from Littlewood \citep[p. 479]{Montgomery:2007} that 
\begin{align*}
	\pi(x) - \operatorname{li}(x) = \Omega_{\pm}\left(x^{1/2}  (\log x)^{-1}\log \log \log x\right).
\end{align*}
The growth is slow, however, and it has been conjectured by Monach and Montgomery~\citep[p. 484]{Montgomery:2007} that 
\begin{align*}
	\pi(x) - \operatorname{li}(x) = O\left(x^{1/2}  (\log x)^{-1}(\log \log \log x)^2\right).
\end{align*}
Numerically, we observe that $\operatorname{Var}[\pi_{\operatorname{RM}}(x)]$ grows slightly faster than 
the latter estimate squared (SI Figure S2), but how fast needs to be investigated further.

\subsubsection{The RM and its relation to the Riemann hypothesis}

The compelling agreement between the RM and the primes is particularly relevant in understanding why the yet unresolved Riemann Hypothesis (RH) should be true. From a classical result by von Koch~\citep{vonKoch:1901ui}, it follows that 
this conjecture can be stated on the form
${|\pi(x) - \operatorname{Ri}(x)|^2 = O(x (\log x)^2)}$, allowing for a direct comparison with \eqref{eq:primevar}.
%
%
%
We know that 
${
	\sum_{i=1}^{[x]} \epsilon_i^2 \sim x /\log x
}$, 
so necessarily, what makes or breaks the RH is the behaviour of the covariance term
 ${
	2 \sum_{i=1}^{[x]} \sum_{i<j}^{[x]} \epsilon_i \epsilon_j
}$.
The prediction stemming from the RM, however, is that this term is always negative, hence ensuring the RH with good margin. As the RM reveals, the negative covariance term--characteristic of subrandom sequences---is an immediate consequence of the recursive structure of the primes. More precisely, it follows from Monach and Montgomery's conjecture above, as well as from our numerical inquiry (SI Figure S2), that the leading term of
${
	2 \sum_{i=1}^{[x]} \sum_{i<j}^{[x]} \epsilon_i \epsilon_j
}$
should be $-x/\log x$. In contrast to this, one should note that proving 
${
	2 \sum_{i=1}^{[x]} \sum_{i<j}^{[x]} \epsilon_i \epsilon_j \leq c x (\log x)^2  
}$
for some constant c and $x$ large enough would be enough to prove the RH.

One important observation we made earlier, following from the fact that the elements of the $\operatorname{RM}_c$ are placed in a tightly constrained subspace of the RM, is that all realisations of $\pi_{\operatorname{RM}_c}(x)$ appear to be strongly correlated (\figref{fig5}). From this observation, a reasonable speculation is that the constraint put on the elements in the $\operatorname{RM}_c$ forces 
all realisations of $\pi_{\operatorname{RM}_c}(x)$, including $\pi(x)$,  to strictly satisfy the asymptotic estimates of 
$\mathbf{E}\left[\pi_{\operatorname{RM}_c}(x)\right]$
and 
$\mathbf{E}\left[ | \pi_{\operatorname{RM}_c}(x)-\operatorname{Ri}(x)|^2\right]$.
In particular, we expect that  
$\pi_{\operatorname{RM}_c}(x) \sim x/ \log x$ 
and 
$ | \pi_{\operatorname{RM}_c}(x) - \operatorname{Ri}(x) |^2 =  O(x/ \log x)$ 
for any outcome of the ${\operatorname{RM}_c}$.

Essentially, our inquiry suggests focusing on the covariance term in \eqref{eq:primevar} as a possible strategy for closing in on the RH. Whether this eventually will drive any progress is hard to judge, specifically considering that the many other avenues sought out to crack open the RH have all encountered unmountable obstacles. In the least, however, the RM appears to provide an intuitive basis for understanding why the RH should be true, and it could presumably serve as a constructive starting point for further theoretical explorations. For example would a proof that all realisations of the $\pi_{\operatorname{RM}_c}(x)$ satisfy $|\pi_{\operatorname{RM}_c}(x) - \operatorname{Ri}(x)|^2 = O(x(\log x)^2)$ also be a proof of the RH. An alternative formulation of this direction is provided by the following equivalent conjecture:
\begin{conjecture}[Equivalent RH] 
Assume any  realisation of the RM such that 
$$
|\pi_{\operatorname{RM}}(x) - \operatorname{Ri}(x)|^2 \neq O(x(\log x)^2).
$$
Then there exists at least one value of $n$ such that ${\prod_{k=0}^{\pi(\sqrt{n})} \rho_k(n+a) >  n}$. In other words, the assumed realisation is not an element in the $\operatorname{RM}_c$.
%
%
%
%
%
\end{conjecture}

The RM predicts negative covariance terms also when we replace $\pi(x)$ by either $\pi(x; d, a)$ or $\pi(x; \mathcal{H})$, counting primes in the arithmetic progression $a+nd$ or prime k-tuples $\mathcal{H}$, respectively. In these cases though, the negative covariance is less pronounced, which can be understood directly from the fact that these counting functions count from a subset of the 'primes' for a given realisation of the RM. This supports and explains the claim made by Richard Brent, that "Twin primes (seem to be) more random than primes" \citep{Brent:2014}. One should note that 
any advancement on the RH from the perspective of the covariance term would most likely lead to progress also in these cases. 

\section{Conclusion}

The absence of a fundamental explanation of the prime numbers and their distribution in the natural numbers ranks as a central conceptual problem in mathematics, one that has proven itself resistant for centuries and that continues to tantalise professional as well as amateur mathematicians. Arguably, our conception of the primes are shaped by their standard number theoretical definition, which provides a top-down perspective---the building blocks (the primes) are defined in terms of what they build (the natural numbers).  
Acknowledging the possibility that a bottom-up perspective could provide a valuable complement to the status quo, we have here examined an equivalent definition by Euclid that allows the primes to be generated recursively from first principles---without reference to the natural numbers. The immediate outcome of this reversed perspective is that it reveals the structure of the primes in the natural numbers and in this sense lays bare the organising principle underlying their distribution. As such, our inquiry brings forth a 
long-sought explanation of the ultimate nature of the primes \citep{QuotesWeb, Luque:2009, Sautoy:2017}.

%
%



%
%
%
%
%
In addition, the proposed perspective permits a random model of primes that accounts for their structure as well as their asymptotic density. The tight linkage between the primes and this model establishes the latter as a constructive tool in developing bottom-up and intuitive understanding of problems related to the distribution of primes---including the Riemann hypothesis and the Hardy-Littlewood k-tuple conjecture---as well as in possibly devising new proof strategies for such problems. The importance of having formulated a model with said properties is aptly reflected in two remarks by Andrew Granville and Terence Tao. Granville, commenting on a result by Maier \citep{Maier:1985}, writes in \citep{Granville:1995} that "Presumably we will remain unable to fully understand the finer details until a model is proposed that adequately accounts for both the sieve of Eratosthenes, and Gauss's density statement." Tao, on the other hand, states in his recent lecture notes on analytic number theory \citep{Tao:web:RandomModels} that  "...we do not have a single unified model for the prime numbers (other than the primes themselves, of course)...many of the models about the primes do not fully take into account the multiplicative structure of primes". As substantiated in this paper, our model satisfies the requirements implied by Granville and Tao, which, in light of Tao's remark, lifts it as a candidate for a single unified model for the primes.  

In condensed form, the gist of this paper can be summed up in a single sentence: \emph{All that underlie the distribution of primes are a unit periodic structure and a recursive rule.} Essentially, this places the set of primes in the same realm as that of another reputed mathematical object, namely the Fibonacci sequence. Just as this sequence flows indefinitely from the seed numbers $F_0=0$ and $F_1=1$ and the recurrence relation $F_n = F_{n-1}+F_{n-2}$, the sequence of primes and all its facets spring from a seed sequence of unit period paired with a recurrence relation anchored in Euclid's definition of primes. 
When viewed with this bottom-up lens, the unruly manifestation of the primes in the natural numbers fades into the background as a mere aftereffect of what the sequence of primes truly is---an  architecture of ever-growing layers of periodic sequences.  
%
%
%
%
%
%
%
%
%
%
%
%
%
As we have brought up earlier, this perspective on the primes admits a geometric interpretation, where the starting point is an unnumbered ruler. The fact that one can build a sound understanding of the distribution of primes even from such a plain and concrete outset, void of abstract notation, accentuates the pedagogical  potential inherent in a first-principles perspective on the primes. And specifically, we recognise in this potential the possibility of inspiring "\dots the reinjection of more or less directly empirical ideas" into the field of number theory, once argued by von Neumann to be a vital condition for conserving the freshness and the vitality of mathematics  \citep{Neumann:1947}.

Ultimately, to fully appreciate the importance of a first-principles perspective on the primes,
it is essential to realise the almost mythological stature the prime numbers have in mathematics, suitably typified here by Zagier~\citep{Zagier:1977}:
\vspace{0.2cm}
\begin{myquote}{1cm }
{\small\textit{
...there is no apparent reason why one number is prime and another not. To the contrary, upon looking at these numbers one has the feeling of being in the presence of one of the inexplicable secrets of creation.
}}
\end{myquote}
\vspace{0.2cm}
\noindent This outlook on the primes, along with similar ones~\citep{QuotesWeb, Sautoy:2017}, tells us that at a fundamental level the  
established understanding of these mathematical objects is inadequate.
%
Nevertheless, this gap in understanding, which has outlasted centuries of advances in number theory, ceases to exist when a first-principles perspective is employed---suggesting that the only 'secret' there ever was to the primes was us looking at them from a skewed angle. In closing, therefore, we have allowed for an elementary resolution to an age-old number theoretical 'mystery'.

\vskip6pt
%
%
%
%
%
%


\small
\bibliography{primereferences}
\bibliographystyle{authordate1}
\vspace{2mm}
\end{document}